\documentclass{amsart}
\usepackage{amsmath, amssymb,epic,graphicx,mathrsfs,enumerate}
\usepackage[all]{xy}
\usepackage{color}
\usepackage{comment}
\usepackage{amscd}

\usepackage{amsthm}
\usepackage{amssymb}
\usepackage{latexsym}
\usepackage{longtable}
\usepackage{epsfig}
\usepackage{amsmath}
\usepackage{hhline}

\DeclareMathOperator{\aut}{Aut}

\DeclareMathOperator{\Hom}{Hom}

\DeclareMathOperator{\core}{Core}

\DeclareMathOperator{\alt}{Alt}
 
\DeclareMathOperator{\der}{Der}

\newtheorem{thm}{Theorem}
\newtheorem{cor}[thm]{Corollary}
 \newtheorem{lemma}[thm]{Lemma}
\newtheorem{prop}[thm]{Proposition} 
 \newtheorem{defn}[thm]{Definition}

\numberwithin{equation}{section}

\renewcommand{\footnote}{\endnote}
\newcommand{\ignore}[1]{}\makeglossary

\begin{document}
	\bibliographystyle{amsplain}
	\title[Genus, thickness and crossing number]{Genus, thickness and crossing number
	of graphs encoding the generating properties of  finite groups}
\author{Cristina Acciarri}
\address{
	Cristina Acciarri\\ Department of Mathematics, University of Brasilia, 70910-900 Bras\'ilia DF, Brazil\\ 
	email:acciarricristina@yahoo.it}
	\author{Andrea Lucchini}
	\address{Andrea Lucchini\\ Universit\`a degli Studi di Padova\\  Dipartimento di Matematica \lq\lq Tullio Levi-Civita\rq\rq\\ Via Trieste 63, 35121 Padova, Italy\\email: lucchini@math.unipd.it}

	\begin{abstract} Assume that $G$ is a finite group and let $a$ and $b$ be non-negative integers. We define an undirected  graph $\Gamma_{a,b}(G)$ whose vertices correspond to the elements of $G^a\cup G^b$ and in which two tuples $(x_1,\dots,x_a)$ and $(y_1,\dots,y_b)$ are adjacent if and only $\langle x_1,\dots,x_a,y_1,\dots,y_b \rangle =G.$  Our aim is to estimate the genus, the thickness and the crossing number
		of the graph $\Gamma_{a,b}(G)$ when $a$ and $b$ are positive integers.
	\end{abstract}
	\maketitle

\section{Introduction}
Generating sets of a finite group may be quite complicated. If a group $G$ is $d$-generated, the question of which sets of $d$ elements of $G$ generate $G$ is nontrivial. The simplest interesting case is when $G$ is 2-generated. One tool developed to study generators of 2-generated finite groups is the generating graph $\Gamma(G)$ of $G;$ this is the graph which has the elements of $G$ as vertices and an edge between two elements $g_1$ and $g_2$ if $G$ is generated by $g_1$ and $g_2.$ Note that the generating graph may be defined for any group, but it only has edges if $G$ is $2$-generated.   A wider family of graphs which encode the generating property of $G$ when $G$ is an arbitrary finite group was introduced and investigated in \cite{ac}. The definition of these graphs is the following. Assume that $G$ is a finite group and let $a$ and $b$ be non-negative integers. We define an undirected  graph $\Gamma_{a,b}(G)$ whose vertices correspond to the elements of $G^a\cup G^b$ and in which two tuples $(x_1,\dots,x_a)$ and $(y_1,\dots,y_b)$ are adjacent if and only $\langle x_1,\dots,x_a,y_1,\dots,y_b \rangle =G.$  
Notice that $\Gamma_{1,1}(G)$ is the generating graph of $G,$ so these graphs can be viewed as a natural generalization of the generating graph. 

\

Let $\Delta$ be a graph. The {\slshape{genus}} $\gamma(\Delta)$ of $\Delta$ is the minimum  integer $g$ such that there exists an embedding of $\Delta$ into the orientable surface $S_g$ of genus $g$ (or in other words the minimum number $g$ of handles which must be added to a sphere so that $\Delta$ can be embedded on the resulting surface).
The {\sl{thickness}} $\theta(\Delta)$ of  $\Delta$ is the minimum number of planar graphs into which the edges of $\Delta$ can be partitioned. The {\slshape{crossing number}} $\rm{cr}(\Delta)$ of $\Delta$ 
is the minimum number of crossings in any simple drawing of $\Delta(G).$ In this paper we investigate genus, thickness and crossing number of the graphs $\Gamma_{a,b}(G)$, when
$1\leq a\leq b$ and $a+b\geq d(G),$ where $d(G)$ is the smallest cardinality of a generating set of $G.$ Notice that the case $a=0$ is not interesting: the graph $\Gamma_{0,b}(G)$ is a star with an internal node corresponding to the empty set and with $\phi_G(b)$ leaves, being $\phi_G(b)$ be the number of the generating $b$-uples of $G$. Our main result is the following:
\begin{thm}\label{stima}
Assume that $G$ is a nontrivial $d$-generated  finite group and that $a, b$ are positive integer with $a+b\geq d.$	
Then
$$\begin{aligned}\gamma(\Gamma_{a,b}(G))&\geq \frac{|G|^b}{6}
\left(\frac{\sqrt{|G|}}{16}-3\right),\\
\theta(\Gamma_{a,b}(G))&\geq \frac{\sqrt{|G|}}{48},\\
{\rm{cr}}(\Gamma_{a,b}(G))&\geq \frac{|G|^{d+\frac{1}{2}}}{29}\left(\frac{1}{ 2^{11}}-\frac{70}{|G|^{3/2}}\right).
\end{aligned}$$
\end{thm}
In order to estimate $\gamma(\Gamma_{a,b}(G)),$ $\theta(\Gamma_{a,b}(G))$ and ${\rm{cr}}(\Gamma_{a,b}(G))$, it is important to obtain a lower bound for the ratio
$e(\Gamma_{a,b}(G))/v(\Gamma_{a,b}(G))$ between the number of edges and the number of vertices of the graph $\Gamma_{a,b}(G).$ We will see in Section \ref{conti}, that this is essentially related to the estimation of the ratio $\phi_G(d)/|G|^{d-1}$ for a $d$-generated finite group. Our main result in this direction is the following.
\begin{thm}\label{main}If $G$ is a $d$-generated finite group, then
	$$\frac{\phi_G(d)}{|G|^{d-1}}\geq \frac{\sqrt{|G|}}{2}.$$
\end{thm}
We think that this is a result of independent interest. For example it implies the following corollary.

\begin{cor}\label{coro}Let $G$ be a finite group and let $d=d(G).$ Denote by $\rho(G)$ the number of elements $g$ in $G$ such that $G=\langle g, x_1,\dots,x_{d-1}\rangle$, for some $x_1,\dots,x_{d-1}\in G.$ We have
	$$\rho(G)\geq \frac{|G|^{1-\frac{1}{2d}}}{2^\frac{1}{d}}$$
\end{cor}

Recall that a graph is said to be embeddable in the plane, or planar, if it can be drawn in the plane so that its edges intersect only at their ends. In \cite{planar} a classification of the 2-generated finite groups with planar generating graph is given. We generalize this result as follows.
\begin{thm}\label{pla}
	Let $G$ be a nontrival finite group  and let $a$ and $b$  be two positive integers with $a\leq b$ and $a+b\geq d(G).$ Then $\Gamma_{a,b}(G)$ is planar if and only if one of the following occurs:
	\begin{enumerate}
		\item $G\in \{C_3, C_4, C_6, C_2 \times C_2, D_3, D_4, Q_8, C_4\times C_2,  D_6\}$
	and $(a,b)=(1,1).$
	\item $G\cong C_2$ and either $a=1$ or $(a,b)=(2,2).$
	\end{enumerate}
\end{thm}

\section{Proof of Theorem \ref{main}}

Let $G$ be a $d$-generated finite group and   let  $\phi_G(d)$  denote  the number of the generating $d$-uples $(g_1,\dots,g_d)\in G^d$ with $\langle g_1,\dots,g_d\rangle=G$. Clearly $P_G(d)=\phi_G(d)/|G|^d$ coincides with the probability that $d$ randomly chosen elements from $G$ generate $G.$ 
 \begin{defn}For a  $d$-generated finite group $G$, set
$$\alpha(G,d):=\frac{\phi_G(d)}{|G|^{d-1}}=P_G(d)|G|.$$ 
\end{defn}



Let $N$ be a normal subgroup of a finite group $G$ and choose $g_1,\dots,g_k\in G$ with the property that $G=\langle g_1,\dots,g_k\rangle N.$ By a result of Gasch{\"u}tz \cite{g1} the cardinality of the set $$\Phi_N(g_1,\dots,g_k)=\{(n_1,\dots,n_k)\in N\mid \langle g_1n_1,\dots,g_kn_k\rangle=G\}$$ does not depend on the choice of $g_1,\dots,g_k.$ Let
$$P_{G,N}(k)=\frac{|\Phi_N(g_1,\dots,g_k)|}{|N|^k}.$$ Notice that if $k\geq d(G/N),$ then
$P_{G,N}(k)=P_G(k)/P_{G/N}(k).$
\begin{defn}Let $N$ be a normal subgroup of  a $d$-generated finite group $G.$ Set 
	$$\alpha(G,N,d):=\frac{\alpha(G,d)}{\alpha(G/N,d))}=P_{G,N}(d)|N|.$$
\end{defn}

\begin{lemma}\label{minab}
Assume that $N$ is a minimal abelian normal subgroup of a $d$-generated finite group $G$. We have $|N|=p^a,$ where $p$ is a prime and $a$ is a positive integer. Let $c$ be the number of complements of $N$ in $G$.
Then $$\alpha(G,N,d)=\frac{p^{d\cdot a}-c}{p^{(d-1)\cdot a}}\geq p^a-p^{a-1}=p^{a-1}(p-1).$$
In particular
\begin{enumerate}
	\item $\alpha(G,N,d)=1$ if and only if $|N|=2$, $N$ has a complement in $G$ 
	and $G/N$ admits $C_2^{d-1}$ as an epimorphic image.
	\item $\alpha(G,N,d)\geq 3/2$ if  $|N|=2,$  $N$ has a complement in $G$ 
	and  $C_2^{d-1}$ is not an epimorphic image of $G/N.$
	\item $\alpha(G,N,d)\geq 2$ in all the remaining cases.
\end{enumerate}
\end{lemma}
\begin{proof}By \cite[Satz 2]{g2}, $P_{G,N}(d)=1-c/p^{d\cdot a},$ hence $\alpha(G,N,d)=\frac{p^{d\cdot a}-c}{p^{(d-1)\cdot a}}$. If $c\neq 0,$ then $c$ is the order of the group $\der(G/N,N)$  of derivations from $G/N$ to $N$; in particular $c$ is a power of $p$. Moreover, since $G$ is $d$-generated, it must be $c<p^{d\cdot a}$ and consequently
$$\alpha(G,N,d)=\frac{p^{d\cdot a}-c}{p^{(d-1)\cdot a}}\geq \frac{p^{d\cdot a}-p^{d\cdot a-1}}{p^{(d-1)\cdot a}}=p^a-p^{a-1}.$$ In particular we can have $\alpha(G,N)<2$ only if $|N|=2$ and $c\neq 0.$
Let $H$ be  a complement of $N$ in $G$ and let $K=H^\prime H^2.$ We have $c=|\der(H,N)|=|\Hom(H/K,N)|.$ Since $G$ is $d$-generated, we have $H/K\cong C_2^t$ with $t<d.$ We have $c=2^t$ and $\alpha(G,N,d)=2-2^{t-d+1}.$
\end{proof}


If a group $G$ acts on a group $A$ via automorphisms, then we say that $A$ is a {$G$-group}. If $G$ does not stabilise any nontrivial proper subgroup of $A$, then $A$ is called an {irreducible} $G$-group. Two $G$-groups $A$ and $B$ are said to be {$G$-isomorphic}, or $A\cong_G B$, if there exists a group isomorphism $\phi: A\rightarrow B$ such that 
$\phi(g(a))=g(\phi(a))$ for all $a\in A, g\in G$.  Following  \cite{JL}, we say that two  $G$-groups $A$ and $B$  are {$G$-equivalent} and we put $A \equiv_G B$, if there are isomorphisms $\phi: A\rightarrow B$ and $\Phi: A\rtimes G \rightarrow B \rtimes G$ such that the following diagram commutes:

\begin{equation*}
\begin{CD}
1@>>>A@>>>A\rtimes G@>>>G@>>>1\\
@. @VV{\phi}V @VV{\Phi}V @|\\
1@>>>B@>>>B\rtimes G@>>>G@>>>1.
\end{CD}
\end{equation*}

\

Note that two $G$\nobreakdash-isomorphic
$G$\nobreakdash-groups are $G$\nobreakdash-equivalent. In the abelian case, the converse is true:
if $A_1$ and $A_2$ are abelian and $G$\nobreakdash-equivalent, then $A_1$
and $A_2$ are also $G$\nobreakdash-isomorphic.
It is known (see for example \cite[Proposition 1.4]{JL}) that two  chief factors $A_1$ and $A_2$ of $G$ are  $G$-equivalent if and only if  either they are  $G$-isomorphic, or there exists a maximal subgroup $M$ of $G$ such that $G/\core_G(M)$ has two minimal normal subgroups, $N_1$ and $N_2$,
$G$-isomorphic to $A_1$ and $A_2$ respectively.
Let $A=X/Y$ be a chief factor of $G$. 
 We say that   $A=X/Y$ is a Frattini chief factor if  $X/Y$ is contained in the Frattini subgroup of $G/Y$; this is equivalent to saying that $A$ is abelian and there is no complement to $A$ in $G$.
The  number of non-Frattini chief factors $G$-equivalent to $A$ in any chief series of $G$  does not depend on the series, and so this number is well-defined: we will denote it by  $\delta_A(G)$.

The following numerical results will be useful.
\begin{lemma}\cite[9.15 p. 54]{feller} Let $n>0$, then
	$$\sqrt{2\pi}\cdot n^{n+\frac{1}{2}}\cdot e^{-n}\cdot e^{\frac{1}{12n+1}}\leq n!\leq  \sqrt{2\pi}\cdot n^{n+\frac{1}{2}}\cdot e^{-n}\cdot e^{\frac{1}{12n}}.$$
\end{lemma}

\begin{cor}\label{corollario}If $t<n,$ then
$$\frac{n!}{(n-t)!} \geq \frac{9}{10}
\frac{n^t}{e^{t}}.$$
\end{cor}

\begin{proof}
$$	\begin{aligned}
\frac{n!}{(n-t)!}&\geq 
\frac{n^{n+\frac{1}{2}}}{e^n}\frac{e^{n-t}}{(n-t)^{n-t+\frac{1}{2}}}\cdot  \frac{e^{\frac{1}{12n+1}}}{e^{\frac{1}{12(n-t)}}}
\geq 
\frac{n^{n+\frac{1}{2}}}{e^n}\frac{e^{n-t}}{(n-t)^{n-t+\frac{1}{2}}}\cdot  \frac{1}{e^{\frac{1}{12}}}\geq\\
&\geq \frac{n^{n+\frac{1}{2}}}{e^n}\frac{e^{n-t}}{(n-t)^{n-t+\frac{1}{2}}}\cdot  \frac{9}{10}\geq 
\frac{9}{10}\frac{n^{n+\frac{1}{2}}}{(n-t)^{n-t+\frac{1}{2}}}\cdot \frac{1}{e^{t}}\geq\\
&\geq 
\frac{9}{10}\frac{n^{n+\frac{1}{2}}}{n^{n-t+\frac{1}{2}}}\cdot \frac{1}{e^{t}}=\frac{9}{10}\frac{n^t}{e^{t}}.\qedhere
\end{aligned}$$
\end{proof}

\begin{prop}\label{nonab}
	Let $G$ be a finite group and let $B$ be a non-abelian chief factor of $G.$ Denote by $t=\delta_G(B)$ the number of factors $G$-equivalent to $B$ in a given chief series of $G$. More precisely
	let $X_1/Y_1, X_2/Y_2,\dots,X_t/Y_t$, with $Y_t\leq X_t\leq \dots \leq Y_1\leq X_1$, be the factors $G$-equivalent to $B$ in a given chief series of $G$. For $1\leq i\leq t,$ let
	$\alpha_i=\alpha(G/Y_i,X_i/Y_i,d)$. We have
	$$\prod_{1\leq i\leq t}\alpha_i \geq 	\frac{9}{10}\left(\frac{53|B|}{90e}\right)^t.$$
\end{prop}
\begin{proof}
	Let  $L=G/C_G(B)$ be the monolithic
	primitive group associated to $B$ and assume $L=\langle l_1,\dots,l_d\rangle.$ 
	Moreover define $\Gamma:=C_{\aut (B)} (L/B)|$, $\gamma=|\Gamma|,$ $\Phi:=\Phi_B(l_1,\dots,l_d).$ By \cite[Proposition 16]{crowns}, for  $1 \leq i \leq t$, we have
	$$\alpha_i=\frac{|\Phi|}{|B|^{d-1}}-\frac{(i-1)\gamma}{|B|^{d-1}}.$$
	Let $\rho=|\Phi|/\gamma$ (notice that $\rho$ is an integer) and let $\tau=|B|^{d-1}/\gamma.$ It follows  from  \cite[Theorem 1.1]{lon} that 
	$\rho/\tau \geq \frac{53}{90}|B|.$
	In view of Corollary \ref{corollario}  we have
	$$\prod_{1\leq i\leq t}\alpha_i = \frac{\rho(\rho-1)\cdots (\rho-(t-1))}{\tau^t}\geq \frac{9}{10\cdot e^{t}}\left(\frac{\rho}{\tau}\right)^t\geq 
		\frac{9}{10}\left(\frac{53|B|}{90e}\right)^t. \qedhere$$
\end{proof}


Next we deal with the proof of Theorem \ref{main}.
\begin{proof}[Proof of Theorem \ref{main}]
	Let $X_t\leq X_{t-1}\leq \dots \leq X_1=G$ be a chief series of $G$ and for $1\leq i\leq t-1,$ let $\alpha_i=\alpha(G/X_{i+1},X_i/X_{i+1},d).$ Since $d(G)=d,$ it must be $\delta_G(C_2)\leq d$ and this implies in particular that
	there exists at most a unique index $j^*$ such that $X_{j^*}/X_{j^*+1}$ has order 2, is complemented in $G/X_{j^*+1}$ and the  quotient $G/X_{j^{*}}$ admits $C_2^{d-1}$ as an epimorphic image. If $|X_i/X_{i+1}|=2$ and $i\neq j^*,$ then, by Lemma \ref{minab}, $\alpha_i\geq {3}/{2}\geq \sqrt 2 = \sqrt{|X_i/X_{i+1}|}.$ If $X_i/X_{i+1}$ is abelian and
	$|X_i/X_{i+1}|=p_i^{n_i}>2,$  then, again  by Lemma \ref{minab}, $\alpha_i\geq p_i^{n_i-1}(p_i-1)\geq \sqrt {p_i^{n_i}}=\sqrt{|X_i/X_{i+1}|}.$ Now assume that $B$ is a non-abelian chief factor of $G$ and let $$I_B=\{1\leq k\leq t-1 \mid X_k/X_{k+1}\equiv_G B\}.$$ By Proposition \ref{nonab}, noticing that  $\delta_B(G)=|I_B|$ and  $|B|\geq 6\sqrt{|B|}$ since $|B|\geq 60,$ we have
	$$\begin{aligned}\prod_{k\in I_B}\alpha_k&\geq \frac{9}{10}\left(\frac{53|B|}{90e}\right)^{\delta_B(G)}
	\geq \left(\frac{53|B|}{100e}\right)^{\delta_B(G)}\geq\\
	&\geq  \left(\frac{|B|}{6}\right)^{\delta_B(G)}\geq \left(\sqrt{|B|}\right)^{\delta_B(G)}=\prod_{k\in I_B}\sqrt{|X_k/X_{k+1}|}. 
	\end{aligned}
	$$
	
	The result follows since $\alpha(G,d)=\prod_{1\leq i \leq t-1}\alpha_i$ and $|G|=\prod_{1\leq i \leq t-1}|X_i/X_{i+1}|$.
\end{proof}

We close this section with the proof of Corollary \ref{coro}.
\begin{proof}[Proof of Corollary \ref{coro}]
	By Theorem \ref{main},
	$$\rho(G)^d\geq \phi_G(d)=\alpha(G,d)|G|^{d-1}\geq \frac {|G|^{\frac{1}{2}}|G|^{d-1}}{2}=
	\frac {|G|^{d-\frac 1 2}}{2}.\qedhere
	$$
\end{proof}

\section{Proof of Theorem \ref{stima}}\label{conti}

Before proving Theorem \ref{stima}, we recall some general results in graph theory concerning lower bounds for the genus, the thickness and the crossing number of a simple graph $\Delta$.
\begin{prop}\cite[7.2.4 -  F35]{hand}\label{P1}
	If  $\Delta$ is a simple graph with $e$ edges and $v$ vertices, then 
	$$\gamma(\Delta)\geq 1-\frac{v}{2}+\frac{e}{6}\geq \frac{v}{6}\left(\frac{e}{v}-3\right).$$
\end{prop}
\begin{prop}\cite[10.3.6 (a)]{bm}\label{P2}.
	If  $\Delta$ is a  simple graph with $e$ edges and $v\geq 3$ vertices, then 
$$\theta(\Delta)\geq \frac{e}{3v-6}.$$
\end{prop}
\begin{prop}\cite[Theorem 6]{ack}\label{P3}
	If  $\Delta$ is a  simple graph with $e$ edges and $v$ vertice, then 
$${\rm{cr}}(\Delta)\geq \frac{e^3}{29v^2}-\frac{35}{29}v.$$
\end{prop}

 Assume that $G$ is a finite group and let $a$ and $b$ be positive integers. 
 Let $d=a+b\geq d(G).$ If $a\neq b$ then $\Gamma_{a,b}(G)$ is a bipartite graphs with two parts, one corresponding to the elements of $G^a$ and the other to the elements of $G^b$. In particular $\Gamma_{a,b}(G)$ has $|G|^a+|G|^b$ vertices and there exists a bijective correspondence between  the set of the generating $d$-uples of $G$ and the set of the edges of $\Gamma_{a,b}(G)$: indeed if $\langle g_1,\dots,g_d\rangle=G,$ then $(g_1,\dots,g_a)$ and $(g_{a+1},\dots,g_{d})$ are adjacent  vertices of the graph. Hence the number of edges of $\Gamma_{a,b}(G)$  is  $\phi_G(d).$ The situation is different if $a=b.$ In that case $\Gamma_{a,a}(G)$ has $|G|^a$ vertices,
 $\phi_G(a)$ loops and other $(\phi_G(d)-\phi_G(a))/2$ edges connecting two different vertices (in other words if $e$ is the the number of edges, excluding the loops, and $l$ is the number of loops, then $2e+l=\phi_G(d)$);
 indeed the two elements $(g_1,\dots,g_a,g_{a+1},\dots,g_d)$ and 
 $(g_{a+1},\dots,g_d,g_1,\dots,g_a)$ give rise to the same edge in $\Gamma_{a,a}(G)$.
 Summarizing, let $\nu$ and $\eta$ be, respectively, the number of vertices and edges of $\Gamma_{a,b}(G),$ excluding the loops. We have
$$|G|^b\leq \nu\leq |G|^a+|G|^b\leq 2|G|^{d-1}.$$ Moreover 
$\eta=\phi_{G}(a+b)$ if $a\neq b$, $\eta=(\phi_{G}(2a)-\phi_G(a))/2$ if $a=b.$ If $\phi_G(a)\neq 0,$ then $\phi_{G}(2a)\geq \phi_G(a)|G|^a$, so $\phi_G(a)\leq \phi_{G}(2a)/|G|^a.$ So if $|G|\geq 2,$ then $\eta\geq \phi_{G}(d)/4.$
By applying Theorem \ref{main} and  Propositions \ref{P1},\ref{P2}  and \ref{P3}  respectively it follows that if $G\neq 1,$ then we have the following inequalities.
$$\gamma(\Gamma_{a,b}(G))\geq \frac{\nu}{6}\left(\frac{\eta}{\nu}-3\right)\geq \frac{|G|^b}{6}
\left(\frac{\phi_G(d)}{8|G|^{d-1}}-3\right)\geq \frac{|G|^b}{6}
\left(\frac{\sqrt{|G|}}{16}-3\right).
$$
$$\theta(\Gamma_{a,b}(G))\geq \frac{\eta}{3\nu}\geq
\frac{\phi_G(d)}{24|G|^{d-1}}\geq \frac{\sqrt{|G|}}{48}.
$$
$$\begin{aligned}{\rm{cr}}(\Gamma_{a,b}(G))&\geq \frac{\eta^3}{29\cdot \nu^2}-\frac{35}{29}\cdot \nu\geq \frac{(\phi_G(d))^3}{29\cdot 4^3 \cdot 4 \cdot (|G|^{d-1})^2}-\frac{70\cdot|G|^{d-1}}{29}\\
&\geq \frac{\phi_G(d)|G|}{29\cdot 4^5}-\frac{70\cdot|G|^{d-1}}{29}\geq \frac{|G|^{d+\frac{1}{2}}}{29\cdot 2^{11}}-\frac{70\cdot|G|^{d-1}}{29}.
\end{aligned}$$
This concludes the proof of Theorem \ref{stima}.

\section{Proof of Theorem \ref{pla}}

The main goal of this section is to prove Theorem \ref{pla}. We star with two preliminary results.

\begin{prop}\label{crit}\cite[Lemma 9.23]{pa}.
	A  simple bipartite planar graph on $v$ vertices, whose every connected component contains at least three vertices, can have not more than $2v-4$ edges.
\end{prop}	

\begin{lemma}\label{tre}
	Let $G$ be a finite group and let $b\geq d(G).$ Consider the set $W=\{(x_1,\dots,x_b)\in G^b\mid \langle x_1,\dots,x_b\rangle=G\}.$ If $G$ is not cyclic, then $|W|\geq 3.$
\end{lemma}
\begin{proof}
	Assume $d=d(G)$ and $G=\langle g_1,\dots,g_d\rangle.$ Then $(g_1,g_2,g_3,\dots,g_d,1,\dots,1),$ $(g_1g_2,g_2,g_3,\dots,g_d,1,\dots,1)$ and $(g_1,g_1g_2,g_3,\dots,g_d,1,\dots,1)$ are three different elements of $W.$
\end{proof}

We are now ready to  embark on  the proof of  Theorem \ref{pla}.

\begin{proof}[Proof of Theorem \ref{pla}] Let $a$ and $b$ be positive integers with $a+b\geq d(G).$ We want to discuss when $\Gamma_{a,b}(G)$ is planar.
	We assume $a+b\geq d(G)$ and $a\leq b.$ If $a=0$, then $\Gamma_{a,b}(G)$ is a star, so it is planar. We may exclude from our discussion the case $a=b=1$, since the result in this case follows from the main result in \cite{planar} (notice that the cyclic group $C_5$ appears in the statement of \cite[Theorem 1.1]{planar} but not in the statement of Theorem \ref{pla}: this is because in \cite{planar} the identity element is not included in the vertex-set of $\Gamma_{1,1}(G)).$

	First assume that $G=\langle g \rangle$ is cyclic. 
	\begin{itemize}
		\item If $a\geq 3,$ take
		$$\begin{aligned}\alpha_1=(1,1,g,1,\dots,1), \alpha_2=(1,g,g,1,\dots,1), \alpha_3=(1,g,1,1,\dots,1) \in G^a,\\
		\beta_1=(g,1,g,1,\dots,1), \beta_2=(g,g,g,1,\dots,1), \beta_3=(g,g,1,1,\dots,1) \in G^b.
		\end{aligned}$$ 
		\item 	If $a=2$ and $|G|\neq 2$, take
		$$\begin{aligned}\alpha_1&=(1,g), \alpha_2=(g,1), \alpha_3=(g,g) \in G^2,\\
		\beta_1&=(1,g^2,1,\dots,1), \beta_2=(g^2,1,\dots,1), \beta_3=(g^2,g^2,1,\dots,1) \in G^b.
		\end{aligned}$$ 
		\item If $a=2$ and $|G|=2$ and $b\geq 3,$ take
		$$\begin{aligned}\alpha_1&=(1,g), \alpha_2=(g,1), \alpha_3=(g,g) \in G^2,\\ 
		\beta_1&=(1,g,g,1\dots,1), \beta_2=(g,1,g,1\dots,1), \beta_3=(g,g,g,1\dots,1) \in G^b.
		\end{aligned}$$
	\end{itemize}
	In all these cases, since $\alpha_i$ and $\beta_j$ are adjacent for every $1\leq i, j \leq 3,$ $\Gamma_{a,b}(G)$ contains $K_{3,3},$ so it is not planar.
	If $a=b=2$ and $|G|=2,$ then $\Gamma_{2,2}(G)\cong K_4$ is planar.
	If $a=1$ and $|G|>2,$ then we may consider the subgraph of $\Gamma_{1,b}(G)$ induced by the following vertices:
	$(1),(g),(g^2), (g,x,\dots,x)\in G^b$ for $x\in G.$ This subgraph is bipartite with $3+|G|$ vertices and $3|G|$ egdes. Since $3|G|>2(3+|G|)-4,$ it follows from Proposition \ref{crit}, that this graph is not planar. On the other hand,  if $a=1$ and $|G|=2,$ then it can be easily seen that the graph $\Gamma_{1,b}(G)$ is planar.
	
	\
	
	Now assume that $G$ is not cyclic. Let $d=d(G)$ and $G=\langle g_1,\dots,g_d\rangle.$
	
	\
	
	\noindent First assume that $a\geq 2.$
	If $a+b=d$, then set
	$$\begin{aligned}
	\alpha_1&=(g_1,g_2,g_3,\dots,g_a)\in G^a,\\ \alpha_2&=(g_1,g_1g_2,g_3,\dots,g_a)\in G^a,\\  \alpha_3&=(g_1g_2,g_2,g_3,\dots,g_a)\in G^a,\\
	\beta_1&=(g_{a+1},g_{a+2},g_{a+3},\dots,g_b)\in G^b,\\ \beta_2&=(g_{a+1}g_{a+2},g_{a+2},g_{a+3},\dots,g_b)\in G^b,\\ \beta_3&=(g_{a+1},g_{a+1}g_{a+2},g_{a+3},\dots,g_b)\in G^b.
	\end{aligned}$$
	If $a+b>d$, choose three different elements $x,y,z$ of $G$ and set
	$$\begin{aligned}
	\alpha_1&=(g_1,g_2,g_3,\dots,g_a)\in G^a,\\ \alpha_2&=(g_1,g_1g_2,g_3,\dots,g_a)\in G^a,\\  \alpha_3&=(g_1g_2,g_2,g_3,\dots,g_a)\in G^a,\\
	\beta_1&=(g_{a+1},g_{a+2},g_{a+3},\dots,g_b,x,\dots,x)\in G^b,\\ \beta_2&=(g_{a+1},g_{a+2},g_{a+3},\dots,g_b,y,\dots,y)\in G^b,\\ \beta_3&=(g_{a+1},g_{a+2},g_{a+3},\dots,g_b,z,\dots,z)\in G^b.
	\end{aligned}$$
	In both cases, since $\alpha_i$ and $\beta_j$ are adjacent for every $1\leq i, j \leq 3,$ $\Gamma_{a,b}(G)$ contains $K_{3,3},$ so it is not planar.
	
	\			
	
	\noindent Assume $a=1$ and $a+b>d.$
	Let $W=\{(x_1,\dots,x_b)\in G^b\mid \langle x_1,\dots,x_b\rangle=G\}$ and let $x,y,z$ be three different elements of $G.$ We may consider the subgraph  of  $\Gamma_{1,b}(G)$ induced by following vertices:
	$(x),(y),(z), w\in W$. This subgraph  is bipartite with $3+|W|$ vertices and $3|W|$ egdes. Since, by Lemma \ref{tre},  $|W|\geq 3,$ it follows $3|W|>2(3+|W|)-4,$ and consequently, by Proposition \ref{crit}, this graph is not planar.
	
	\
	
	\noindent Finally assume $a=1$ and $a+b=d.$ Let $H=\langle g_2,\dots,g_b\rangle.$ If $H$ is cyclic, then $d(G)\leq 2,$ in contradiction with $d(G)=1+b$ and $b>1.$
	Let $x,y,z$  be three different elements of $H$ and
	$W=\{(x_1,\dots,x_b)\in G^b\mid \langle x_1,\dots,x_b\rangle=H\}$. We may consider the subgraph  of  $\Gamma_{1,b}(G)$ induced by following vertices:
	$(g_1x),(g_1y),(g_1z), w\in W$. It is bipartite with $3+|W|$ vertices and $3|W|$ egdes. Since $H$ is not cyclic, we have $|W|\geq 3$ by Lemma \ref{tre}. 
	It follows $3|W|>2(3+|W|)-4,$ and consequently,  by Proposition \ref{crit}, this graph is not planar.
\end{proof}


\begin{thebibliography}{99} 
	
	
	\bibitem{ac} C. Acciarri and A. Lucchini, Graphs encoding the generating properties of a finite group, Mathematische Nachrichten, to appear. arXiv:1707.08348.
	
	\bibitem{ack} E. Ackerman
	On topological graphs with at most four crossings per edge,
	Comput. Geom. 85 (2019), 101574, 31 pp.
	
	
	\bibitem{bm} J. A. Bondy and U. S. R. Murty, {{Graph theory}}, Graduate Texts in Mathematics, 244. Springer, New York, 2008. 
	
	
	
	
	
	
	
	
	
	
	
	
	
	\bibitem{pa} M. Cygan, F. Fomin, L.  Kowalik, D. Lokshtanov, D. Marx, M. Pilipczuk,  M. Pilipczuk, S. Saurabh, Parameterized algorithms. Springer, Cham, 2015 xviii+613
	
	\bibitem{crowns} E. Detomi and A. Lucchini, Crowns and factorization of
	the probabilistic zeta function of a finite group,  {{J. Algebra}} {265} (2003), no. 2, 651--668.
	
	\bibitem{lon} E. Detomi and A. Lucchini,  Probabilistic generation of finite groups with a unique minimal normal subgroup, {{J. Lond. Math. Soc.}} (2) 87 (2013), no. 3, 689--706.
	
	
\bibitem{feller} W. Feller,  An introduction to probability theory and its applications. Vol. I. Third edition John Wiley \& Sons, Inc., New York-London-Sydney 1968 xviii+509 pp.	
	
	
	
	\bibitem{g1} {W. Gasch\"utz}, Zu einem von B.H. und H. Neumann gestellten Problem, {{Mathematische Nachrichten}} 14 (1955), 249--252.
	
	\bibitem{g2}
	W. Gasch{\"u}tz, {Die {E}ulersche {F}unktion endlicher
		aufl\"osbarer
		{G}ruppen}, {{Illinois J. Math.}} {3} (1959), 469--476. 
	
	\bibitem{hand} Handbook of graph theory. Second edition. Edited by J. Gross, J. Yellen and P. Zhang. Discrete Mathematics and its Applications (Boca Raton). CRC Press, Boca Raton, FL, 2014. xx+1610
		
	\bibitem{JL} P. Jim\'{e}nez-Seral and J. Lafuente, On complemented non-abelian chief factors of a finite group. Israel J. Math. {106} (1998), 177--188. 
	
	
	
	
	
	
	


\bibitem{planar} A. Lucchini, Finite groups with planar generating graph, Australas. J. Combin. 76 (2020), part 1, 220--225.
	
	
	
\end{thebibliography}
\end{document}